\documentclass{amsart}
\usepackage{amsmath,amsfonts,amscd,amssymb}
\usepackage{ lamsarrow, pb-diagram, pb-lams,   latexsym, ulem}
\usepackage{amsmath,amsthm,amssymb}
\newtheorem{theorem}{Theorem}[section]

\newtheorem{lemma}[theorem]{Lemma}
\newtheorem{sublemma}[theorem]{Sublemma}
\newtheorem{corollary}[theorem]{Corollary}
\newtheorem{proposition}[theorem]{Proposition}

\theoremstyle{definition}

\theoremstyle{remark}
\newtheorem{remark}[theorem]{Remark}

\newcommand{\lem}[2]{\begin{lemma}\label{#1} #2\end{lemma}}

\newcommand{\prop}[2]{\begin{proposition}\label{#1}#2\end{proposition}}
\newcommand{\thm}[2]{\begin{theorem}\label{#1}#2\end{theorem}}
\newcommand{\cor}[2]{\begin{corollary}\label{#1}#2\end{corollary}}

\newcommand{\pf}[1]{\begin{proof}#1 \end{proof}}
\newcommand{\pfc}[2]{\begin{proof}[Proof of #1]#2\end{proof}}

\newcommand{\cal}{\mathcal}

\newcommand{\al}{\alpha}
\newcommand{\be}{\beta}
\newcommand{\ga}{\gamma}
\newcommand{\Ga}{{\Gamma}}
\newcommand{\de}{\delta}
\newcommand{\De}{{\Delta}}

\newcommand{\ze}{\zeta}
\newcommand{\et}{\eta}
\renewcommand{\th}{\theta}

\newcommand{\ph}{\varphi}

\newcommand{\ps}{\psi}

\newcommand{\Om}{{\Omega}}
\newcommand{\wt}[1]{\widetilde{#1}}
\newcommand{\ovl}[1]{\overline{#1}}

\newcommand{\Sp}{{\cal S}_{(p)}}

\newcommand{\Nn}{\mathbb N}
\newcommand{\Z}{\mathbb Z}
\newcommand{\Zp}{{\mathbb Z}_{(p)}}
\newcommand{\Q}{\mathbb Q}

\newcommand{\lk}{\langle}
\newcommand{\rk}{\rangle}
\newcommand{\cg}{\equiv}
\def\o+{\oplus}
\newcommand\Op{\bigoplus}
\def\+{{\textstyle \bigoplus}}

\newcommand{\x}{\times}
\newcommand{\ox}{\otimes}

\newcommand{\e}{{\rm Ext}}
\newcommand{\lnb}{\refstepcounter{theorem}\leqno(\thetheorem)}
\newcommand{\lnr}[1]{\lnb\label{#1}}
\newcommand{\nr}{\refstepcounter{theorem}\thetheorem}
\newcommand{\kko}[1]{(\ref{#1})}
\newcommand{\mbx}[1]{\quad\mbox{#1}\quad}
\newcommand{\cf}{{\it cf.}\ }
\newcommand{\sus}{\subset }

\newcommand{\ar}[1]{\xrightarrow{#1} }
\newcommand{\arr}{\to}

\newcommand{\Coker}{{\rm Coker}\ }

\newcommand{\barr}{\begin{array}}
\newcommand{\ear}{\end{array}}
\newcommand{\bsk}{\begin{array}{rcl}}
\newcommand{\esk}{\end{array}}
\newcommand{\skh}{\begin{eqnarray*}}
\newcommand{\sko}{\end{eqnarray*}}
\newcommand{\bcss}{\begin{cases}}
\newcommand{\ecss}{\end{cases}}
\newcommand{\sk}[1]{$$\begin{array}{rcl}#1\end{array}$$}
\newcommand{\skr}[2]{$$\begin{array}{rcl}#1\end{array}\lnr{#2}$$}

\newcommand{\AR}[2]{$$\begin{array}{#1}#2\end{array}$$}

\newcommand{\cass}[1]{\begin{cases}#1\ecss}

\newcommand{\Nd}[1]{\!\!\!\begin{diagram}\node{#1}\end{diagram}\!\!\!}

\newcommand{\dT}[1]{\dgTEXTARROWLENGTH = #1em }

\newcommand{\cn}{&\equiv &}

\newcommand{\en}{&= &}
\newcommand{\es}{\ = \ }

\newcommand{\lrk}[1]{\left\lk #1\right\rk}
\newcommand{\qand}{\mbx{and}}
\newcommand{\eR}{\et_R}

\newcommand{\AN}{Adams-Novikov}
\renewcommand{\ss}{spectral sequence}
\newcommand{\ANSS}{Adams-Novikov spectral sequence}

\newcommand{\cln}{\colon}
\renewcommand{\O}[1]{\ovl{#1}}
\newcommand{\LR}[1]{\left(#1\right)}

\newcommand{\ak}{\quad}

\newcommand{\ki}{&&\ak}

\newcommand{\tp}[1]{t_1^{p^{#1}} }
\newcommand{\p}[1]{^{p^{#1}} }

\newcommand{\C}[2]{{#1\choose #2}}

\newcommand{\Fk}[1]{F\lrk{#1}} 

\newcommand{\dtw}[1]{
\ifnum#1=0  t_1\ox t_1^p  \fi 
\ifnum#1=1  t_1^p\ox t_1^{p^2} \fi 
\ifnum#1>1  \count1=1  \advance \count 1 by #1 
t_1\p{#1}\ox t_1\p{\the\count1}\fi} 

\newcommand{\dctw}[1]{
\ifnum#1=0  t_1^p\ox t_1  \fi 
\ifnum#1=1  t_1^{p^2}\ox t_1^p \fi 
\ifnum#1>1  \count1=1  \advance \count 1 by #1 
t_1\p{\the\count1}\ox t_1\p{#1}\fi} 

\newcommand{\dtt}[1]{
\ifnum#1=0  t_1\ox t_2^p+t_2\ox t_1\p2  \fi 
\ifnum#1=1  t_1^p\ox t_2^{p^2}+t_2^p\ox t_1\p3 \fi 
\ifnum#1>1  \count1=1 \count 2=2 \advance \count 1 by #1 \advance \count 2 by #1
t_1\p{#1}\ox t_2\p{\the\count1}+t_2\p{#1}\ox t_1\p{\the\count2}\fi} 

\newcommand{\dtf}[1]{
\ifnum#1=0  t_1\ox t_3^p+t_2\ox t_2\p2+t_3\ox t_1\p3  \fi 
\ifnum#1=1  t_1^p\ox t_3^{p^2}+t_2^p\ox t_2\p3+t_3^p\ox t_1\p4 \fi 
\ifnum#1>1 \count1=1 \count 2=2 \count 3=3 \advance \count 3 by #1 \advance \count 1 by #1 \advance \count 2 by #1
t_1\p{#1}\ox t_3\p{\the\count1}+t_2\p{#1}\ox t_2\p{\the\count2}+t_3\p{#1}\ox t_1\p{\the\count3}\fi} 

\newcommand{\liminj}{\varinjlim}

\makeatletter
\def\varholim@#1#2{%
  \vtop{\m@th\ialign{##\cr
    \hfil$#1\operator@font holim\hfil\cr
    \noalign{\nointerlineskip\kern1.5\ex@}#2\cr
    \noalign{\nointerlineskip\kern-\ex@}\cr}}%
}
\def\holimproj{%
  \mathop{\mathpalette\varholim@{\leftarrowfill@\textstyle}}\nmlimits@
}
\def\holiminj{%
  \mathop{\mathpalette\varholim@{\rightarrowfill@\textstyle}}\nmlimits@
}
\makeatother

\newcommand{\arT}[2]{\ {\dT{#1} \ar{#2}}\ }

\newcommand{\der}{\partial} 
\newcommand{\bm}[1]{\mbox{\boldmath $#1$}
} 

\makeatletter
\newcommand{\xRightarrow}[2][]{%
\ext@arrow 0055{\Rightarrowfill@}{#1}{#2}%
}
\def\Rightarrowfill@{\arrowfill@\Relbar\Relbar\Rightarrow}
\newcommand{\xLeftarrow}[2][]{%
\ext@arrow 0055{\Leftarrowfill@}{#1}{#2}%
}
\def\Leftarrowfill@{\arrowfill@\Leftarrow\Relbar\Relbar}
\newcommand{\xLongleftrightarrow}[2][]{%
\ext@arrow 0055{\llrafill@}{#1}{#2}%
}
\def\llrafill@{\arrowfill@\Leftarrow\Relbar\Rightarrow}
\makeatother

\begin{document}

\title[The first line of the Bockstein spectral sequence]{
The first line of the Bockstein spectral sequence on a monochromatic spectrum at an odd prime
}
\author{ Ryo Kato}
\address{Graduate school of Mathematics, Nagoya University, Aichi, 464-8601, Japan}
\email{ryo\_kato\_1128@yahoo.co.jp}
\author{ Katsumi Shimomura}
\address{Department of Mathematics, Faculty of Science, Kochi University, Kochi, 780-8520, Japan}
\email{katsumi@math.kochi-u.ac.jp}
\date{}
\subjclass{Primary 55T99; Secondary 55Q45}
\maketitle

\begin{abstract}
	The chromatic spectral sequence is introduced in \cite{mrw} to compute the $E_2$-term of the \ANSS\ for computing the stable homotopy groups of spheres.
The $E_1$-term $E_1^{s,t}(k)$ of the spectral sequence is an Ext group of $BP_*BP$-comodules. 
There are a sequence of Ext groups $E_1^{s,t}(n-s)$ for non-negative integers $n$ with $E_1^{s,t}(0)=E_1^{s,t}$, and  Bockstein spectral sequences computing a module $E_1^{s,*}(n-s)$ from $E_1^{s-1,*}(n-s+1)$.
So far, a small number of the $E_1$-terms are determined.
Here, we determine the $E_1^{1,1}(n-1)=\e^1M^1_{n-1}$ for $p>2$ and $n>3$ by computing the Bockstein spectral sequence with $E_1$-term $E_1^{0,s}(n)$ for $s=1,2$.
As an application, we study the non-triviality of the action of $\alpha_1$ and $\beta_1$ in the homotopy groups of the second Smith-Toda spectrum $V(2)$.
\end{abstract}

\section{Introduction}

Let $p$ be a prime number, $\Sp$ the stable homotopy category of $p$-local spectra, and $S$ the sphere spectrum localized at $p$.
Understanding homotopy groups $\pi_*(S)$ of $S$ is one of the principal problems in stable homotopy theory.
The main vehicle for computing $\pi_*(S)$ is the Adams-Novikov spectral sequence based on the Brown-Peterson spectrum $BP$.
$BP$ is the $p$-typical component of $MU$, the complex cobordism spectrum, and that it has homotopy groups $BP_*=\pi_*(BP)=\Z_{(p)}[v_1,v_2,\cdots]$ where $v_n$ is a canonical generator of degree $2p^n-2$.
In order to study the $E_2$-term of the \ANSS, H. Miller, D. Ravenel and S. Wilson \cite{mrw} introduced the chromatic \ss.
It was designed to compute the $E_2$-term, but has the following deeper connotation. 
Let $L_n\cln \Sp \to \Sp$ denote the Bousfield-Ravenel localization functor with respect to $v_n^{-1}BP$ (\cf \cite{r:loc}).
It gives rise the  chromatic filtration $\Sp\to\cdots\to L_n\Sp\to L_{n-1}\Sp\to\cdots\to L_0\Sp$ of the stable homotopy category of spectra, which is a powerful tool for understanding the category. 
The chromatic $n$th layer of the spectrum $S$ can be determined from the homotopy groups of $L_{K(n)}S$, the Bousfield localization of $S$ with respect to the $n$th Morava $K$-theory $K(n)$ that it has homotopy groups  $K(n)_*=v_n^{-1}\Z/p[v_n]$ for $n>0$ and $K(0)_*=\Q$. 
By the chromatic convergence theorem of Hopkins-Ravenel \cite{r:orange}, $S$ is the inverse limit of the $L_nS$.
Let $E(n)$ be 
the $n$th Johnson-Wilson spectrum $E(n)$ with $E(n)_*=v_n^{-1}\Zp[v_1,\cdots,v_n]$ for $n>0$ and $E(0)=K(0)$. 
It is Boufield equivalent to $v_n^{-1}BP$ and also to $K(0)\vee\cdots\vee K(n)$,  i.e. $L_{E(n)}=L_n=L_{K(0)\vee\cdots\vee K(n)}$.
We notice that
$E(0)=H{\mathbb Q}$, the rational Eilenberg-MacLane spectrum, and $E(1)$ is the $p$-local Adams summand of  periodic complex $K$-theory.
Futhermore, $E(2)$ is closely  related to elliptic cohomology.
So far, we have no geometric interpretation of homology theories $K(n)$ or $E(n)$ when $n>2$.

From now on, we assume that the prime $p$ is odd.
We explain the $E_1$-term of the chromatic spectral sequence.
The Brown-Peterson spectrum $BP$ is a ring spectrum that induces the Hopf algebroid
$
	(BP_*, BP_*(BP))=(BP_*, BP_*[t_1,$ $t_2,\dots])
$
in the standard way \cite{r:book}, 
and we have an induced Hopf algebroid
$$
	(E(n)_*,E(n)_*(E(n)))=(E(n)_*, E(n)_*\ox_{BP_*} BP_*(BP)\ox_{BP_*}E(n)_*)
$$
where $E(n)_*$ is considered to be a $BP_*$-module by sending $v_k$ to zero for $k>n$.
Then, the $E_1$-term is given by
$$
	E_1^{s,t}(n-s)=\e_{E(n)_*(E(n))}^t(E(n)_*, M^s_{n-s}).
$$
Here, 
$M^s_{n-s}$ denotes the $E(n)_*(E(n))$-comodule $E(n)_*/(I_{n-s}+(v_{n-s}^\infty,v_{n-s+1}^\infty,\dots,$ $v_{n-1}^\infty))$, in which
$I_k$ denotes the ideal of $E(n)_*$ generated by $v_i$ for $0\le i<k$ ($v_0=p$), and  $M/(w^\infty)$ for $w\in E(n)_*$ and an $E(n)_*$-module $M$ denotes the cokernel of the localization map $M\to w^{-1}M$.
In order to study the stable homotopy groups $\pi_*(L_{K(n)}S)$,
we study here the homotopy groups of the monochromatic component $M_nS$ of $S$ (see \cite{r:loc}).
Then, the $E_2$-term $E_2^{s,t}(M_nS)$ of the \ANSS\ for computing $\pi_*(M_nS)$ is the $E_1$-term $E_1^{n,s}(0)$ of the chromatic \ss.
In \cite{mrw}, the authors also introduced the $v_{n-s}$-Bockstein spectral sequence 
$E_1^{s-1,t+1}(n-s+1)\Rightarrow E_1^{s,t}(n-s)$ associated to a short exact sequence 
$$0\to M^{s-1}_{n-s+1}\ar{\ph} M^s_{n-s}\ar{v_{n-s}} M^s_{n-s}\to 0$$ 
of $E(n)_*(E(n))$-comodules, where $\ph(x)=x/v_{n-s}$.
So far, the $E_1$-term $E_1^{s,t}(n-s)$ is determined in the following cases (\cf \cite{r:book}):
\sk{
	(s,t,n)\en (0,t,n) \mbx{for (a) $n\le 2$, (b) $n=3$, $p>3$, (c) $t\le 2$  by Ravenel \cite{r:M}, }\\
		\ki \mbx{(Henn \cite{h} for $n=2$ and $p=3$),}\\
		\en (1,0,n) \mbx{for  $n\ge 0$ by Miller, Ravenel and Wilson \cite{mrw},}\\
		\en (s,t,n) \mbx{for  $n\le 2$ by Shimomura and his colaborators: Arita \cite{as},}\\
	\ki  \mbx{ Tamura \cite{st}, Yabe \cite{sy} and Wang \cite{sw},
	(\cite{s:86}, \cite{s:97}, \cite{s:00}),}\\
		\en (1,1,3) \mbx{by Shimomura \cite{s-tottori}, Hirata and Shimomura \cite{hs},}\\
		\en (2,0,n) \mbx{for $n>3$ by Shimomura \cite{s:90}, for $n=3$ by Nakai \cite{n:n}, \cite{n:k}.}
}
In this paper, we determine the structure of  $E_1^{1,1}(n-1)$ for $n>3$. 
The case $n=3$, which is special, is treated in \cite{s-tottori} and \cite{hs}.
The result is the first step to understand $\pi_*(L_{K(n)}S)$ for $n>3$ as explained above.
We proceed to state the result.

In this paper, we consider only the cases $s=0$ and $s=1$, and,
 hereafter, put
$$
	v=v_n\qand u=v_{n-1}.
$$
Furthermore, we put
$$
	F=\Z/p,
$$
and consider the coefficient ring $K(n)_*=F[v_n^{\pm 1}]=F[v^{\pm 1}]=E(n)_*/I_n$, 
\sk{
	A\en E(n)_*/I_{n-1}\qand 
	B\es M_{n-1}^1\es  A/(u^\infty)\es \Coker (A\to u^{-1}A).
}

\noindent
Since the ideal $I_{n-1}$ is invariant, $(A, \Ga)=(A, E(n)_*(E(n))/I_{n-1})$ is a Hopf algebroid, and we use the abbreviation
$$
	\e^sM=\e^s_\Ga(A,M)
$$
for a $\Ga$-comodule $M$.
Then, the chromatic $E_1$-terms are
$$
	E_1^{0,t}(n)=\e^tK(n)_*\qand E_1^{1,t}(n-1)=\e^tB.
$$
We have the $u$-Bockstein spectral sequence 
$$
	E_1=\e^*K(n)_*\Longrightarrow \e^*B
\lnr{bss}$$
associated to the short exact sequence
$$
	0\arr K(n)_* \ar{\ph} B \ar{u} B \arr 0,
\lnr{bss}
$$
where $\ph$ is a homomorphism defined by $\ph(x)=x/u$.

Let $R$ be a ring, and let $R\lrk{g}$ denote the $R$-module generated by $g$.
The $E_1$-term of the $u$-Bockstein spectral sequence was  determined by Ravenel \cite{r:M} as follows:
\thm{E1}{ $\e^0K(n)_*=K(n)_*$ and 
\sk{
	\e^1K(n)_*\en K(n)_*\lrk{h_i, \ze_n: 0\le i <n},\\
	\e^2K(n)_*\en K(n)_*\lrk{\ze_nh_i, b_{i}, g_i, k_i, h_jh_k: 0\le i <n, \ 0\le j< k-1<n-1}.
}}

In the theorem, the generators $h_i$ and $b_i$ are represented by $\tp i$ and $\sum_{k=1}^{p-1}\frac1p\C p k t_1^{kp^i}\ox t_1^{(p-k)p^i}$ of the cobar complex $\Om_\Ga^*K(n)_*$, respectively, and $g_i$ and $k_i$ are given by the Massey products 
$$
	g_i=\lrk{h_i,h_i,h_{i+1}}\qand k_i=\lrk{h_i,h_{i+1},h_{i+1}}.
\lnr{Massey}$$
In order to determine the module $\e^0B$, Miller, Ravenel and Wilson \cite{mrw} introduced elements $x_{i}$ and integers $a_{i}$ 
in \cite[(5.11) and (5.13)]{mrw}, where they denoted them by $x_{n,i}$ and $a_{n,i}$,  such that $x_{i}\cg v^{p^i}$ mod $I_n$ with the action of the connecting homomorphism $\de$ given in \cite[(5.18)]{mrw}:
$$
	\de(v^s/u)= sv^{s-1}h_{n-1} \qand
	\de(x_{i}^s/u^{a_{i}})= sv^{(sp-1)p^{i-1}}h_{[i-1]} \ \mbox{for $i\ge 1$.}
\lnr{de:1}
$$
Hereafter, we let
$$
[i]\in\{0,1,\dots, n-2\}
$$
be the principal representative of the integer $i$ module $n-1$.
The elements $x_i$ and the integers $a_i$ are defined inductively by $x_0=v$ and $a_0=1$, and for $i>0$,
\skr{
	x_i\en \cass{x_{i-1}^p&\mbox{for $i= 1$ or $[i]\ne 1$,}\\ x_{i-1}^p-u^{b_{n,i}}v^{p^i-p^{i-1}+1}&\mbox{for $i>1$ and $[i]=1$, and}}\\
	a_{i}\en \cass{pa_{i-1}&\mbox{for $i= 1$ or $[i]\ne 1$,}\\
	pa_{i-1}+p-1&\mbox{for $i>1$ and $[i]=1$.}
}
}{xi,ai}
Here, $b_{n,k(n-1)+1}=(p^n-1)(p^{k(n-1)}-1)/(p^{n-1}-1)$. 
The result \kko{de:1} determines the differentials of the Bockstein spectral sequence, which
implies:

\thm{H0M1n}{{\rm (\cite[Th. 5.10]{mrw})} As a $k_*$-module, 
$$
	\e^0B=L_\infty\o+ \Op_{p\nmid s, i\ge 0} L_{a_{i}}\lrk{x_{i}^s}.
$$
Here, $	k_*=k(n-1)_*=F[u]$, $L_i=k_*/(u^i)$  
 and
$
 L_\infty=k_*/(u^\infty)=\liminj_i L_i
$.
}

This theorem together with \kko{de:1} 
implies the following:

\cor{cor:coker}{The cokernel of $\de\colon \e^0B\to \e^1K(n)_*$ is the $F$-module generated by
\sk{
	v^t\ze_n,&&
	v^{tp-1}h_{n-1}, \ak h_j\mbx{for $0\le j<n-1$, and}\\
	v^{sp^k}h_j&& \mbox{for $0\le j<n-1$, where $[k]\ne [j]$, $s\not\cg -1$ $(p)$, or $s\cg -1$ $(p^2)$,}
}
for integers $s$  and $t$ with $p\nmid s$.}

By Theorem \ref{E1},
the module $\e^1K(n)_*$ is the direct sum of $\ze_n\e^0K(n)_*=\ze_nK(n)_*$, $F\lrk{ h_j}$ for $j\in \Z/(n-1)$ and the modules 
$$
	V_{(i,j,s)}=F\lrk{v^{sp^i}h_j}
$$ for $(i,j,s)\in\Nn \x \Z/n\x \O \Z$.
Here, $\Nn$ denotes the set of non-negative integers, and $\O \Z=\Z\setminus p\Z$.
We partition $\Nn \x \Z/n$ as follows:

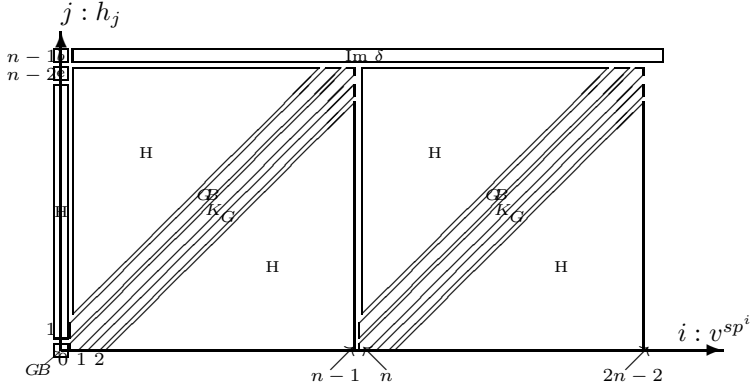
\setlength\unitlength{.8truecm}
\begin{picture}(10,7)(-0.5,-1)
\thicklines
	\put(0,0){\vector(1,0){11}}
  \put(0,0){\vector(0,1){5.3}}
  \put(-0.05,-0.25){\scriptsize $0$}
  \put(0.25,-0.25){\scriptsize $1$}
  \put(0.55,-0.25){\scriptsize $2$}
  \put(4.45,-0.25){$\nearrow$}
  \put(4.15,-0.5){\scriptsize $n-1$}
  \put(5,-0.25){$\nwarrow$}
  \put(5.3,-0.5){\scriptsize $n$}
  \put(9,-0.5){\scriptsize $2n-2$}
  \put(9.35,-0.25){$\nearrow$}

  \put(-0.25,0.25){\scriptsize $1$}
  \put(-.9,4.8){\scriptsize $n-1$}
  \put(-.9,4.5){\scriptsize $n-2$}
  \put(0,5.5){$j:h_j$}
  \put(10.25,0.15){$i:v^{sp^i}$}
\thinlines
  \put(-.1,4.8){\framebox(.2,.2){\tiny $b$}}
  \put(.2,4.8){\framebox(9.8,.2){\tiny Im $\de$ }}
  \put(-.1,4.5){\framebox(.2,.2){\tiny e}}
  \put(-.1,.2){\framebox(.2,4.2){\tiny H}}
  \put(-.1,-.1){\framebox(.2,.2){}}
  \put(-.6,-.4){\tiny $G\!B$}
  \put(-.3,-.2){\tiny $\nearrow$}

\put(.2,4.7){\line(1,0){4.1}}
\put(.2,.6){\line(0,1){4.1}}
\put(.2,.6){\line(1,1){4.1}}
\put(1.3,3.2){\tiny H}

\put(.15,.2){\line(0,1){.25}}
\put(.16,.2){\line(1,1){4.5}}
\put(.16,.45){\line(1,1){4.25}}
\put(4.4,4.7){\line(1,0){.25}}
\put(2.25,2.5){\tiny $G\!\!B$}

\put(.3,0){\line(1,1){4.58}}
\put(.15,.1){\line(1,1){4.6}}
\put(4.75,4.7){\line(1,0){.13}}
\put(.15,0){\line(0,1){.12}}
\put(4.88,4.57){\line(0,1){.13}}
\put(2.4,2.25){\tiny $K$}

\put(.45,0){\line(1,1){4.42}}
\put(.65,.0){\line(1,1){4.22}}
\put(4.88,4.225){\line(0,1){.2}}
\put(2.65,2.13){\tiny $G$}

\put(.75,0){\line(1,1){4.13}}
\put(4.88,0){\line(0,1){4.13}}
\put(3.4,1.3){\tiny H}

\put(4.8,0){\begin{picture}(4,4)(0,0)

\put(.2,4.7){\line(1,0){4.1}}
\put(.2,.6){\line(0,1){4.1}}
\put(.2,.6){\line(1,1){4.1}}
\put(1.3,3.2){\tiny H}

\put(.15,.2){\line(0,1){.25}}
\put(.16,.2){\line(1,1){4.5}}
\put(.16,.45){\line(1,1){4.25}}
\put(4.4,4.7){\line(1,0){.25}}
\put(2.25,2.5){\tiny $G\!\!B$}

\put(.3,0){\line(1,1){4.58}}
\put(.15,.1){\line(1,1){4.6}}
\put(4.75,4.7){\line(1,0){.13}}
\put(.15,0){\line(0,1){.12}}
\put(4.88,4.57){\line(0,1){.13}}
\put(2.4,2.25){\tiny $K$}

\put(.45,0){\line(1,1){4.42}}
\put(.65,.0){\line(1,1){4.22}}
\put(4.88,4.225){\line(0,1){.2}}
\put(2.65,2.13){\tiny $G$}

\put(.75,0){\line(1,1){4.13}}
\put(4.88,0){\line(0,1){4.13}}
\put(3.4,1.3){\tiny H}
\end{picture}}
\end{picture}

\noindent
More precisely,
\sk{
	H\en \{(0,j): 1\le j<n-2\}\\
	\ki  \cup \{(i,j): i>0,\ [i]\ne n-3, n-2,\ 2+[i]\le j\le n-2\}\\
	\ki \cup \{(i,j): i>0,\ [i]\ne 0, 1,\ 0\le j\le [i]-2\},\\
	GB\en \{(i,[i]): i\ge 0\},\\
	K\en \{(i,[i]-1): i>0, \ [i]\ne 0\}\qand \\
	G\en \{(i,[i]-2): i>1, \ [i]\ne 0,1\}.
}
We  introduce notation
\sk{
	V_{(0,n-2)}\en \Op_{s\in \O\Z'}V_{(0,n-2,s)},\\
	V_{(0,n-1)}\en \Op_{t\in\Z} V_{(0,n-1,tp-1)}\es F[v^{\pm p}]\lrk{v^{-1}h_{n-1}},\\
	C_X\en \Op_{(i,j)\in X,\ s\in\O\Z}V_{(i,j,s)} \mbx{for a subset $X\sus \Nn\x \Z/n$,}\\
	\O C_{GB}\en \Op_{(i,j)\in GB}\LR{\LR{\Op_{s\in \O{\O\Z}}V_{(i,j,s)}}\o+ \LR{\Op_{t\in\Z} V_{(i,j,tp^2-1)}}}\\
	\en \Op_{(i,[i],s)\in \wt{GB}}V_{(i,j,s)}\o+ \Op_{i\ge 0}F[v^{\pm p^{i+2}}]\lrk{v^{-p^{i}}h_{[i]}}\qand \\
	C_O\en F\lrk{\th, h_j:j\in\Z/(n-1)}.}
Here, for $e(i)\es (p^i-1)/(p-1)$, $\th=v^{e(n-2)}h_{n-2}$, 
\sk{
	\O \Z'\en \O \Z\setminus\{e(n-2)\}, \ak 
	\O{\O \Z}\es \{n\in \O\Z:p\nmid (s+1)\} \qand \\
	\wt{GB}\en \{(i,[i],s): s\in \O{\O\Z}\}.
}
We also consider the subset $\bm T$ of $\Nn\x\Z/n\x \O \Z$ defined by
\sk{
	\bm T\en \{ (i,j,s)\in \Nn\x\Z/n\x\O \Z: \mbox{ $p\nmid (s+1)$ or $p^2\mid (s+1)$ if $[i]=j$,}\\
	\ki \qquad \mbox{$p\mid (s+1)$ if $(i,j)=(0,n-1)$, and $s\ne e(n-2)$ if $(i,j)=(0,n-2)$}\}.
}
In this notation, the cokernel of $\de$ in Corollary \ref{cor:coker} is given by
\skr{
	\Coker \de\en \ze_nK(n)_*\oplus C_O\o+ \Op_{(i,j,s)\in \bm T}V_{(i,j,s)}\\
	&\hspace{-.5in}=& \hspace{-.3in}\ze_nK(n)_*\oplus C_O\o+ V_{(0,n-2)}\o+ V_{(0,n-1)}\o+ C_H\o+ C_K\o+ C_G\o+ \O C_{GB}\\
}{coker}
Finally, we consider the $k_*$-modules:
\sk{
	W_{(i,j,s)}\en L_{a(i,j,s)}\lrk{x_i^sh_j},\\
	W_{(0,n-2)}\en \Op_{s\in \O\Z'}W_{(0,n-2,s)},\\
	W_{(0,n-1)}\en \Op_{t\in\Z} W_{(0,n-1,tp-1)},\\
	B_X\en \Op_{(i,j)\in X,\ s\in\O\Z}W_{(i,j,s)} \mbx{for a subset $X\sus \Nn\x \Z/n$, }\\
	\O B_{GB}\en \Op_{(i,j)\in GB}\LR{\LR{\Op_{s\in \O{\O\Z}}W_{(i,j,s)}}\o+ \LR{\Op_{t\in\Z} W_{(i,j,tp^2-1)}}} \qand\\ 
	C_\infty\en \LR{K(n-1)_*/k_*}\lrk{\th, h_j:j\in\Z/(n-1)}.
}
 Here, $a(i,j,s)$ denotes an integer defined as follows:
for $(i,j)=(0,n-2)$, $a{(0,n-2,s)}=2$ if $p\nmid s(s-1)$, and
\sk{
	a{(0,n-2,s)}\en \cass{
		a_l &  p\nmid t,\ l>0 ,\ [l]\ne 0, n-2,\\
		a_l+e(n-2)+p^{n-3} &  p\nmid t,\ l>0,\ [l]= n-2,\\
		a_l+1 &  p\nmid t,\ l>0,\ [l]= 0}
}
if $s=tp^{l}+e(n-2)$; 
for $(i,j)\in \{(0,n-1)\}\cup H\cup K\cup G\cup GB$, 

\sk{
	a{(i,j,s)}\en \cass{
		p-1 &(i,j)=(0,n-1),\\
		a_i&(i,j)\in H,\\
		a_i+a_{i-1} &(i,j)\in K\cup G,\\
		2a_i&(i,j,s)\in \wt{GB},\\
		(p-1)a_{i+1}& (i,j)\in GB,\ p^2\mid (s+1).\\
		}
}

\thm{main}{The chromatic $E_1$-term $\e^1B=\e^1M^1_{n-1}$ is canonically isomorphic to the $k_*$-module 
$$
	\ze_n\e^0B\o+ C_\infty\o+ W_{(0,n-2)}\o+  W_{(0,n-1)}\o+ B_H\o+ B_K\o+ B_G\o+ \O B_{GB}.
$$
}



Let $V(n)$ be the $n$th Smith-Toda spectrum defined by $BP_*(V(n))=BP_*/I_{n+1}$.
As an application of the theorem, we study the action of $\al_1$ and $\be_1$ on the elements $u^t$ $(t>0)$ in the \AN\ $E_2$-term $E_2^*(V(n))$ in section 6.
In particular, it leads us an geometric result for $n=4$.
In \cite{t}, Toda constructed the self map $\ga$ on $V(2)$ to show the existence of $V(3)$ for the prime $p>5$.
We notice that $\ga^ti\in \pi_*(V(2))$ for the inclusion $i\cln S\to V(2)$ to the bottom cell is detected by $u^t=v_3^t\in BP_*(V(2))$ in the \ANSS.

\thm{1}
{Let $p>5$.
Then $\gamma^ti\alpha_1$ and $\gamma^ti\beta_1$ are nontrivial in $\pi_*(V(2))$ for $t>0$.
}

\section{Bockstein spectral sequence}

We compute the Bockstein spectral sequence by use of the following lemma.

\lem{key}{Let $\de\cln \e^sB\to \e^{s+1}K(n)_*$ be the connecting homomorphism associated to the short exact sequence \kko{bss}. Suppose that $\Coker \de=\Op_k V_k\sus \e^1K(n)_*$ and $\Op_k U_k\sus \e^2K(n)_*$ for $F$-modules $V_k$ and $U_k$,  and there exist $u$-torsion $k_*$-modules $W_k$ fitting in a commutative diagram
\[
\begin{CD}
0@>>> V_k @>{\varphi _*'}>> W_k @>{u}>> W_k @>{{\de}'}>>U_k
\\
@. @VVV @V{f_k}VV @VV{f_k}V @VVV
\\
0@>>> {\rm Coker}\ \de @>{\varphi _*}>> {\e}^1B @>{u}>> {\e}^1B @>{\de}>>\e ^2K(n)_*
\end{CD}
\]
of exact sequences.
Then, $\e^1B=\Op_k W_k$.}

This follows immediately from \cite[Remark 3.11]{mrw}.

Let $\wt\th$ be an element of Corollary \ref{m:3}.
Then, $\wt\th/u^k$ and $h_j/u^k$ for $j\in \Z/(n-1)$ belong to $\e^1B$, and
we define the map $f\cln C_\infty \to \e^1B$ by $f((u^{-k})\th)=\wt\th/u^k$ and $f((u^{-k})h_j)=h_j/u^k$ for $(u^{-k})\in K(n-1)_*/k_*$, so that 
 the short exact sequence
$$
	0\arr C_O\arT2{1/u} C_\infty \ar{u} C_\infty \arr 0
\lnr{infty}$$
 yields a summand of Lemma \ref{key}.

Note that if a cocycle $z$ represents $\ze_n$, then so does $z^p$.
Therefore, we have $\ze_n/u^j\in \e^1B$ represented by $z^{p^j}/u^j$.
The exact sequence \kko{bss} induces the exact sequence
$
	0\arr \e^0K(n)_*\ar{\ph_*} \e^0B\ar{u}\e^0B\ar{\de}\e^1K(n)_*
$, and we have an exact sequence
$$
	0\arr \ze_n\e^0K(n)_*\ar{\ph_*} \ze_n\e^0B\ar{u}\ze_n\e^0B\ar{\de}\ze_n\e^1K(n)_*,
\lnr{sum:1}$$
which is a summand of Lemma \ref{key}.
Together with \kko{infty} and \kko{sum:1},
 Theorem \ref{main} follows from Lemma \ref{key} if  the following sequence is exact for each $(i,j,s)\in \bm T$:
$$
	0\arr V_{(i,j,s)}\ar{\ph_*'} W_{(i,j,s)} \ar{u} W_{(i,j,s)}\ar{\de'} U_{(i,j,s)},
\lnr{ses:ijs}$$
where $U_{(i,j,s)}$ denotes an $F$-module generated by a single generator as follows:
for $(i,j)=(0,n-2)$, $U_{(0,n-2,s)}=\Fk{v^{s-2}k_{n-2}}$ if $p\nmid s(s-1)$,
\sk{
	U_{(0,n-2,s)}\en \cass{
		\Fk{v^{s-p^{l-1}}h_{[l-1]}h_{n-2}} &  p\nmid t, \ l>0 ,\ [l]\ne 0, n-2,\\
		\Fk{v^{s-p^{l-1}}b_{2n-5}} &  p\nmid t, \ l>0,\ [l]= n-2,\\
		\Fk{v^{s-p^{l-1}-1}g_{n-2}} &  p\nmid t, \ l>0,\ [l]= 0;
}
}
if $s=tp^{l}+e(n-2)$;
for $(i,j)\in \{(0,n-1)\}\cup H\cup K\cup G\cup GB$,
\sk{	U_{(i,j,s)}\en \cass{
		\Fk{v^{s-p+1}b_{n-1}} &(i,j)=(0,n-1),\\
		F\lrk{v^{(sp-1)p^{i-1}}h_{[i-1]}h_j}&(i,j)\in H,\\
		\Fk{v^{(s-2)p}k_{n-1}}&(i,j)=(1,0)\in K,\\
		\Fk{v^{(sp^2-p-1)p^{i-2}}k_{[i-2]}}&(i,j)\in K,\ i>1,\\
		\Fk{v^{(sp^2-p-1)p^{i-2}}g_{[i-2]}}&(i,j)\in G,\\
		\Fk{v^{s-p-1}g_{n-1}}&(i,j,s)\in \wt{GB},\ i=0,\\
		\Fk{v^{(sp-2)p^{i-1}}g_{[i-1]}}&(i,j,s)\in \wt{GB},\ i> 0,\\
		F\lrk{v^{(s+1-p)p^i}b_j}& (i,j)\in GB, \ p^2\mid (s+1).
		}
}

Since the mapping $\bm T \to \{U_{(i,j,s)}: (i,j,s)\in \bm T\}$ assigning $(i,j,s)$ to $U_{(i,j,s)}$ is an injection, we see the following:

\lem{inc}{The direct sum of $\ze_n\e^1K(n)_*$ and  $U_{(i,j,s)}$ for $(i,j,s)\in \bm T$ is a sub-$F$-module of $\e^2K(n)_*$. }

The homomorphism $f_k$ in Lemma \ref{key} on $W_{(i,j,s)}$ for $(i,j,s)\in\bm T$ is explicitly given by
$$
	f_{(i,j,s)}(x)=x/u^{a{(i,j,s)}}.
$$
It follows that the homomorphism $\de'$ on it is given by the composite $\de (1/u^{a{(i,j,s)}})$.
Hereafter we denote it by $\de'_{(i,j,s)}$, that is, $\de'_{(i,j,s)}=\de (1/u^{a{(i,j,s)}})$, and
 consider a condition:\\

\noindent
$(\nr\label{cond})_{(i,j,s)}$
	\qquad $\de_{(i,j,s)}'(x)=y$ for the generators $x\in W_{(i,j,s)}$ and $y\in U_{(i,j,s)}$.\\

Note that $\ph'_*(\O x)=u^{a{(i,j,s)}-1}x$ for the generators $\O x\in V_{(i,j,s)}$ and $x\in W_{(i,j,s)}$, since $f_k\ph'_*(\O x)=\ph_*(\O x)=x/u$. Then,

\lem{key:2}{For each $(i,j,s)\in \bm T$, if the condition $\kko{cond}_{(i,j,s)}$ holds, then $\kko{ses:ijs}$  for $(i,j,s)$ is exact and yields a summand of Lemma \ref{key}.}

The relations in \kko{de:1} show immediately\\[-2mm]

\noindent
(\nr\label{sum:2})
{\it The condition $\kko{cond}_{(i,j,s)}$ holds for $(i,j)\in H$.
}

\pfc{Theorem \ref{main}}{The theorem follows from Lemmas \ref{key}, \ref{inc} and \ref{key:2} together with \kko{infty}, \kko{sum:1}, \kko{sum:2},  Lemmas \ref{GB:1}, \ref{GB:2}, \ref{GK} and \ref{0,n-2}, in which the lemmas are proved below.
Indeed, the direct sum of $\ze_n\e^0K(n)_*$, $C_O$ and $V_{(i,j,s)}$ for $(i,j,s)\in\bm T$ is the cokernel of $\de$ by \kko{coker}.  
}

\section{The summands on  $V_{(0,n-1)}$ and $\O C_{GB}$}

We begin with stating some formulae on the Hopf algebroid $(A,\Ga)$:
\skr{
	0\en vt_k\p n +ut_{k+1}\p{n-1}-u\p{k+1}t_{k+1}-t_k\eR(v\p k)\in \Ga\mbx{for $k<n$,}\\
	\eR(u)\en u,\ak \eR(v)\es v+u\tp{n-1}-u^pt_1,\\
	\De(t_k)\en \sum_{i=0}^k t_i\ox t_{k-i}\p i \mbx{for $k<n$, and}\\
	\De(t_n)\en \sum_{i=0}^n t_i\ox t_{n-i}\p i-ub_{n-2}.
}{formulae}

Then the connecting homomorphism $\de\cln \e^1B\to \e^2K(n)_*$ is computed by the differential $d\cln \Om_\Ga^1A\to \Om_\Ga^2A$ of the cobar complex modulo an ideal, which is defined by
$$
	d(x)=1\ox x-\De(x)+x\ox 1.
\lnr{d}$$
We also use the differential $d\cln \Om_\Ga^0A\to \Om_\Ga^1A$ defined by $d(w)=\eR(w)-\et_L(w)$.
For $w,w'\in \Om^0_{\Ga}A$ and $ x \in \Om^1_{\Ga}A$,
these differentials satisfy
\skr{
	d(ww')\en d(w)\eR(w')+wd(w'),\
	 d(wx)\es d(w)\ox x+wd(x),\ {\rm and}\\
	&& \hspace{-.6in}d(x\eR(w))\es d(x)\eR(w)-x\ox d(w).
}{d:rel}

We also use the Steenrod operations $P^0$ and $\be P^0$ on $\e^* C(j)$ for $j\ge 1$ and $\e^* B$ (\cf \cite{may:st}, \cite{r:book}).
Here, $C(j)$ denotes the comodule $A/(u^j)$, and we notice that $C(1)=K(n)_*$.
Let $\wt \Om^s M=\Om^s_{E(n)_*(E(n))}M$ for an $E(n)_*(E(n))$-comodule $M$.
Given a cocycle $x(j)$ of $\wt \Om^s C(j)$, 
$\wt x(j)$  denotes a cochain of $\wt \Om^sE(n)_*$ such that  $\pi_j(\wt x(j))=x(j)$ for the projection $\pi_j\cln \wt \Om^sE(n)_*\to \wt \Om^sC(j)$.
Since $x(j)$ is a cocycle, $d(\wt x(j)^p)=py_j+\sum_{i=1}^{n-2}v_i^pz_{j,i}+u^{jp}z_{j,n-1}$ for some elements $y_j$ and $z_{j,i}\in \wt \Om^{s+1}E(n)_*$.
Under this situation, the Steenrod operations are
 defined by 
\AR{c}{
	P^0([x(j)])=[x(j)^p]\qand \be P^0([x(j)])= [y_{j}] \in \e^{*}C(jp), \qand\\
	P^0([x(j)/u^j])=[x(j)^p/u^{jp}]\qand \be P^0([x(j)/u^j])= [y_{j}/u^{jp}] \in \e^{*}B.
}
Here, $[x]$ denotes the homology class represented by a cocycle $x$.
In particular,
the operation
acts on our elements as follows:
\skr{
	 \be P^0(x_{i}/u^{a_i})\en \cass{v^{p-1}h_{n-1}/u^{p-1}&i=0,\\ 
	x_{i-1}^{p^{2}-1}h_{[i-1]}/u^{(p-1)a_i}&i>0,} \mbx{in $\e^1B$;}
}{Ph0}
\vspace{-.15in}
\skr{	P^0(x_{i}^sh_k/u^j)\en \cass{x_{i+1}^sh_{k+1}/u^{jp}&k\ne n-2,\\ x_{i+1}^sh_{0}/u^{jp-p+1}&k=n-2;}\mbx{in $\e^1B$; and}\\
	 \be P^0(x_{i}^sh_k)\en x_{i+1}^sb_{k}\mbx{in $\e^2K(n)_*$.}
}{Ph}
The following is a folklore (cf. \cite[Corollary A1.5.5]{r:book}):
$$
	 P^0\de =\de P^0 \qand \be P^0\de=-\de\be P^0 \mbx{in $\e^*K(n)_*$.}
\lnr{P0-de}
$$

\lem{GB:1}{The condition $\kko{cond}_{(i,j,s)}$ holds for each $(i,j,s)\in \{(0,n-1,tp-1), (i,j,tp^2-1):t\in \Z, (i,j)\in GB\}$.} 

\pf{For $k\ge -1$, consider a generator $x(k,t)=x_{k}^{tp^{2}-1}h_{[k]}$ for $k\ge 0$ and $x(-1)=x_{0}^{tp-1}h_{n-1}$, and  $\O{(k,t)}$ denotes a triple $(k,[k],tp^2-1)$ if $k\ge 0$ and $(0,n-1,tp-1)$ if $k=-1$.
Then, $(1/u^{ a\O{(k,t)}})(x(k,t))=x_{k+2}^{t-1}\be P^0(x_{k+1}/u^{a_{k+1}})$ for $k\ge -1$ by \kko{Ph0}.
Now, $\de'_{\O{(k,t)}}(x(k,t))$ equals 
$$
  x_{k+2}^{t-1}\de (\be P^0(x_{k+1}/u^{a_{k+1}}))\es - x_{k+2}^{t-1}(\be P^0(x_{k}^{p-1}h_{\O{[k]}}))\es  - x_{k+1}^{\nu (t)}b_{\O{[k]}}
$$
by \kko{P0-de}, \kko{de:1} and \kko{Ph}.
Here, $(\nu(t),\O{[k]})=(tp-1,[k])$ if $k\ge0$ and $=((t-1)p,n-1)$ if $k=-1$.
}

\lem{GB:2}{The condition $\kko{cond}_{(i,[i],s)}$ holds for  $(i,[i],s)\in \wt{GB}$.
}

\pf{We prove this by induction on $i$.
By \kko{formulae} and \kko{d}, we compute mod $(u^3)$
\sk{
	d(v^{s+1-p}t_1\p n)\cn (s+1)uv^{s-p}t_1\p{n-1}\ox \tp n+\C{s+1}2u^2v^{s-p-1}t_1^{2p^{n-1}}\ox \tp n\\
	d((s+1)uv^{s-p}t_2\p{n-1})\cn s(s+1)u^2v^{s-p-1}\tp{n-1}\ox t_2\p{n-1}-(s+1)uv^{s-p}t_1\p{n-1}\ox \tp n
}
to obtain $\de(v^sh_0/u^2)=s(s+1)v^{s-p-1}g_{n-1}$ and so
$$
	\de_{(0,0,s)}'(v^sh_0)\es s(s+1)v^{s-p-1}g_{n-1}.
$$
Apply $P^0$ to it, and we obtain
\sk{	\de_{(1,1,s)}'(v^{sp}h_1)
\en
\de (P^0(v^sh_0/u^2))
\es
 P^0\de (v^sh_0/u^2)
\es
s(s+1) P^0(v^{s-p-1}g_{n-1})\\
\en
s(s+1)v^{sp-p^2-p}g_{n}
\es
s(s+1)v^{sp-2}g_{0}.
}
Here, we notice that $g_n=v^{p^2+p-2}g_0$ in $\e^2K(n)_*$ by \kko{formulae}.
Suppose inductively that $\de'_{(i,1,s)}(x_i^{s}h_1)=s(s+1)v^{(sp-2)p^{i-1}}g_{0}$ for $[i]= 1$, which is \kko{cond}$_{(i,1,s)}$.  
Note that $a_{i+j}=pa_{i+j-1}$ if $0<j<n-2$, and we see that $P^0\de_{(i,j,s)}'=\de'_{(i+1,j+1,s)}P^0$ by
 \kko{P0-de}.
Therefore, $(P^0)^j$ for $j<n-2$ yields the equation for $\de'_{a(i+j,j+1,s)}(x_{i+j}^{s}h_{j+1})$.
At $i'=i+n-2$, for $t=(i',0,s)$,  $\de_{t}'(x_{i'}^{s}h_0)=\de P^0(x_{i'-1}^sh_{n-2}/u^{a(i'-1,n-2,s)})$ (by \kko{Ph}) $=s(s+1) v^{(sp-2)p^{i+n-3}}g_{n-2}$ by \kko{P0-de} and inductive hypothesis.

Note that $a_{i+n-1}=p^{n-1}a_i+p-1$.
Consider the connecting homomorphism $\de_j\cln \e^1M^1_{n-1}\to \e^2C(j)$ associated to the short exact sequence $0\arr C(j)\ar{1/u^j}M^1_{n-1}\ar{u^j} \Nd{M^1_{n-1}}\arr 0$.
Then, $u^{j-1}\de=\de_ju^{j-1}$. Besides, $\de_j(P^0)^k=(P^0)^k\de$ if $p^k\ge j$.
Now  in $\e^2C(p^2+p-1)$, $u^{p^2+p-2}
\de'_{(i+n-1,1,s)}(x_{i+n-1}^{s}h_1)$ equals
\sk{
	&&\hspace{-.3in} u^{p^2+p-2}\de (x_{i+n-1}^{s}h_1/u^{p^{n-1}a+2(p-1)} )
	\es  \de_{p^2+p-1} (P^0)^{n-1}(x_{i}^{s}h_1/u^{a} )\\
	\en (P^0)^{n-1}(s(s+1)v^{(sp-2)p^{i-1}}g_{0})
	\es s(s+1)v^{(sp-2)p^{i+n-2}}g_{n-1}}
for $a=a(i,[i],s)$, which equals
$s(s+1)u^{p^2+p-2}v^{(sp-2)p^{i+n-2}}g_{0}
$ by the relation $u^{p+2}g_{n-1}=u^{p^2+2p}g_0$.
This relation follows from \kko{Massey} and $uh_{n-1}=u^ph_0$ given by $d(v)$.
}

\section{The summands $C_G$ and $C_K$}

We study the action of the connecting homomorphism $\de$
by use of the Massey product. We notice that this is also shown by use of $P^0$-operation considered in the previous section, but we use the Massey product for the sake of simplicity.

\lem{GK}{The condition $\kko{cond}_{(i,j,s)}$ holds for $(i,j)\in G\cup K$.
}

\pf{
We consider the element $(1/u^{a(i,j,s)})(x_i^sh_j)$ the Massey product $\Big\lk sx_{i-1}^{sp-1}/u^{a_{i-1}},$ $ h_{[i-1]}, h_j\Big\rk$.
Then, $\de'_{(i,j,s)}(x_i^sh_j)=
	 \de \lrk{sx_{i-1}^{sp-1}/u^{a_{i-1}}, h_{[i-1]}, h_j}
	= \Big\lk s\de(x_{i-1}^{sp-1}/u^{a_{i-1}}),$ $ h_{[i-1]}, h_j\Big\rk
$, 
 which equals $-\lrk{sv^{sp-2}h_{n-1},h_{0},h_0}=-sv^{(s-2)p}k_{n-1}$ if $i=1$, and $-\lrk{sv^{(sp^2-p-1)p^{i-2}}h_{[i-2]},h_{[i-1]},h_j}=\cass{-sv^{(sp^2-p-1)p^{i-2}}k_{j-1}&j=[i-1],\\-2sv^{(sp^2-p-1)p^{i-2}}g_j& j=[i-2]}$ otherwise.
Here, we note that $\lrk{h_i,h_{i+1},h_i}=2g_i$.
}

\section{The summand $V_{(0,n-2)}$}

Consider
 the elements $c_i=u^{p^{i}}h_{n-1+i}$ and $c_i'=u^{p^{i+1}}h_i$ of $\e^1A$.
The elements have internal degrees $|c_i|=|c_i'|=p^ie(n)q$  
 for $q=2p-2$, and satisfy
$$
	c_i=c_i',\ak c_ic_{i+1}=0,\ak h_{n+i}c_i=0 \qand h_{i+1}c_i=h_{i+1}c_i'=0.
$$

We consider the cochains $\O w_k=u^{e(k-1)}ct_k\p{n-1}$ of the cobar complex $\Om_\Ga^1A$. Then,
\skr{
	\O w_k\en -\O w_{k-1}^p\eR(v)+u^{pe(k-2)}v\p{k-1}ct_{k-1}+u^{p^k+pe(k-2)}ct_k
}{Owk}
for $k>1$ by \kko{formulae}.
Let $w_k$  be a cochain of the cobar complex $\Om_\Ga^1A$ defined inductively by:
\skr{
	w_1\en \tp{n-1}-u^{p-1}t_1\es -\O w_1+u^{p-1} ct_1\qand \\
	w_k\en w_{k-1}^p\eR(v)+(-1)^ku^{pe(k-2)}v\p{k-1}ct_{k-1}
}{wk}
and put
\skr{
	m'_k\en -\sum_{i=1}^{k-1}(-1)^iu\p{i-1}w_{k-i}\p{i}\ox \O w_{i}\qand \\
	m_k\en u\p{k-1}w_k+\sum_{i=1}^{k-1}(-1)^iu\p{i-1}v^{p^ie(k-i)}\O w_i.
}{mk}

\lem{dve}{ $d(v^{e(k)})=m_{k}$. Besides, $d(w_k)=m'_{k}$ if $k\le n$.}

\pf{
We prove the lemma inductively.
Since $d(v)=uw_1=m_1$, we see the case for $k=1$. Indeed, $m_{1}'=0$.

Suppose that the equalities hold for $k-1$. Then, we compute by \kko{d:rel}, \kko{Owk} and \kko{wk},
\sk{
	d(v^{e(k)})\en  d(v^{pe(k-1)})\eR(v)+v^{pe(k-1)}d(v)\\
	&\hspace{-1.2in}=&\hspace{-.7in}  \LR{u\p {k-1}w_{k-1}^p+\sum_{i=1}^{k-2}(-1)^iu\p{i}v^{p^{i+1}e(k-1-i)}\O w_i^p}\eR(v)-uv^{pe(k-1)}\LR{\O w_1-u^{p-1}ct_1}\\
	&\hspace{-1.2in}=&\hspace{-.7in}  u\p {k-1}\LR{w_k-(-1)^ku^{pe(k-2)}v\p{k-1}ct_{k-1}}-uv^{pe(k-1)}\LR{\O w_1-u^{p-1}ct_1}\\
	&\hspace{-1.2in}&\hspace{-.7in}   +\sum_{i=1}^{k-2}(-1)^iu\p{i}v^{p^{i+1}e(k-1-i)}\LR{-\O w_{i+1}+\LR{u^{pe(i-1)}v\p{i}ct_{i}+u^{p^{i+1}+pe(i-1)}ct_{i+1}}},
}
which equals $m_k$,
and similarly,
\sk{
	d(w_k)\en -\sum_{i=1}^{k-2}(-1)^iu\p{i}w_{k-1-i}\p{i+1}\ox \O w_i^p\eR(v)+uw_{k-1}^p\ox (\O w_1-u^{p-1}ct_1)\\
	\ki +(-1)^ku^{pe(k-2)}\LR{u\p{k-1}w_1\p{k-1}\ox ct_{k-1}+v\p{k-1}d(ct_{k-1})}\\
	\en -\sum_{i=1}^{k-2}(-1)^iu\p{i}w_{k-1-i}\p{i+1}\ox \LR{-\O w_{i+1}+\uline{u^{pe(i-1)}v\p{i}ct_{i}+u^{p^{i+1}+pe(i-1)}ct_{i+1}}}\\
	\ki \ +uw_{k-1}^p\ox (\O w_1-\uline{u^{p-1}ct_1} )\\
	\ki +\uline{(-1)^ku^{e(k-2)}\LR{u\p{k-1}w_1\p{k-1}\ox ct_{k-1}+v\p{k-1}d(ct_{k-1})}}
	\es m_k'
}
Here, the underlined terms cancel each other if $k\le n$ by \kko{wk} and \kko{formulae} with the relation $\De (cx)=T(c\ox c)\De (x)$ for the switching map $T\cln \Ga\ox \Ga\to \Ga\ox \Ga$.
}
We also introduce an element
$$
	\O c_k=h_{n+k-1}-u^{(p-1)p^{k}}h_k\in \e^1A.
$$

\cor{cor:dve}{
For each $0<k<n$, the Massey products $\mu_k=\lk u^{p^k},\O c_k, c_{k-1},c_{k-2},$ $\dots,c_1,c_0\rk$ and $\mu_k'=\lrk{\O c_k, c_{k-1},c_{k-2},\dots,c_1,c_0}$ are defined.
In fact, the cocycles
$m_{k+1}$ and $m_{k+1}'$ represent elements of the Massey products $\mu_k$ and $\mu'_k$, respectively.}

In particular, we have

\cor{m:1}{The Massey product $\lrk{u\p{n-3}, \O c_{n-3}, c_{n-4}, \dots, c_0}\sus \e^1A$ is defined and contains zero.}

\lem{m:2}{The Massey product $
	\lrk{\O c_{n-3}, c_{n-4}, \dots, c_0,h_{n-2}}\sus \e^2A
$
contains zero.
}

\pf{
The Massey product 
$
	\lrk{\O c_{n-3}, c_{n-4}, \dots, c_0,h_{n-2}}
$
contains
$$
	\lrk{h_{2n-4}, c_{n-4}, \dots, c_0,h_{n-2}}-\lrk{u^{p^{n-2}-p^{n-3}}h_{n-3}, c_{n-4}, \dots, c_0,h_{n-2}}.
$$
It suffices show that the second term contains zero.
Indeed, the first term does since a defining system cobounds $u^{e(n-3)}ct_{n-1}\p{n-2}$.
Since every Massey product $\lk h_j,h_{j-1},$ $\dots,h_{i+1}, h_i\rk$ for $j-i\leq n-2$ contains zero, all lower products contains zero, and
we see that $\xi=\lrk{h_{n-3},c_{n-4},\dots,c_1,c_0,h_{n-2}}$ is defined.

The statement of \cite[Th.10]{k} itself is applied to our case and says that  there are elements $x_k\in \lrk{c_k,c_{k-1},\dots,c_0,h_{n-2},h_{n-3},c_{n-4},\dots, c_{k+1}}$ for $0\le k\le n-4$, $x_{n-3}\in \lrk{h_{n-3},c_{n-4},\dots,c_1,c_0,h_{n-2}}$ and $x_{n-2}\in \lrk{h_{n-2}, h_{n-3},c_{n-4},\dots,c_1,c_0}$ such that
$
	\sum_{k=0}^{n-2}\pm x_k=0
$.
Its proof tells us that we may take  the elements $x_k$ arbitrary, and
we take $x_k$ so that $x_k=0$ for $0\le k\le n-4$ and $x_{n-2}=0$, whose relations follow from $d(ct_{n-1})$.
Therefore, $x_{n-3}=0$ and the lemma follows.
}

\cor{m:3}{The Massey product $\mu=\lrk{u\p{n-3},\O c_{n-3}, c_{n-4},\dots, c_0,h_{n-2}}$ is defined and contain an element whose leading term is $v^{e(n-2)}h_{n-2}$.}

\lem{0,n-2}{The condition $\kko{cond}_{(i,j,s)}$ holds for $(i,j)=(0,n-2)$. 
}

\pf{
If $p\nmid s(s-1)$, it follows from the computation
\sk{
	d(v^s\tp{n-2})\cn suv^{s-1}\tp{n-1}\ox \tp{n-2}+\C s2 u^2t_1^{2p^{n-1}}\ox \tp{n-2}\mod (u^3)\\
	d(suv^{s-1}ct_2\p{n-2})\cn s(s-1)u^2\tp{n-1}\ox ct_2\p{n-2}-suv^{s-1}\tp{n-1}\ox \tp{n-2} \mod (u^3).
}
Suppose $s=tp^l+e(n-2)$ with $p\nmid t$ and $l>0$.
Let $\wt \th$ denote an element of Corollary \ref{m:3}. 
We take a generator corresponding to $v^{s}h_{n-2}$ to be $v^{s-e(n-2)}\wt\th$.
We denote a representative of $\wt\th$ by $m$, which is congruent to $v^{e(n-2)}t_1\p{n-2}+uv^{pe(n-3)}ct_2\p{n-2}$ modulo $(u^2)$.
Then, $d(v^{s-e(n-2)}m)=tu^{a_l}v^{s-e(n-2)-p^{l-1}}\tp{[l-1]}\ox m\cg tu^{a_l}v^{s-p^{l-1}}\tp{[l-1]}\ox \tp{n-2}$.
This shows the case for $[l]\ne 0, n-2$.

For $[l]= 0$, the similar computation shows that $d(v^{s-e(n-2)}m)\cg tu^{a_l}v^{s-p^{l-1}}(\tp{n-2}\ox \tp{n-2}+uv^{-1}t_1^{p^{n-1}+p^{n-2}}\ox \tp{n-2}+uv^{-1}\tp{n-2}\ox ct_2\p{n-2})$, 
which yields $v^{s-1-p^{l-1}}g_{n-2}$. 
For $[l]=n-2$, $\wt\th h_{n-3}\in u^{e(n-2)}\lrk{h_{2n-4},h_{2n-5},\dots,h_{n-2},h_{n-3}}=\{ u^{e(n-2)+p^{n-3}}b_{2n-5}\}$ in  $C(p^{n-2})$.
Indeed, $u^{e(n-3)}t_n\p{n-3}$ yields the equality by \kko{formulae}.
}

\section{On the action of $\al_1$ and $\be_1$ on Greek letter elements}

In this section, let \ $H^*M$ \ for a  $BP_*(BP)$-comodule \ $M$ \ denote an Ext group $\e^*_{BP_*(BP)}(BP_*,M)$.
Consider the comodule $N_{k-1}(j)=BP_*/(I_{k-1}+(v_{k-1}^j))$ $(v_0=p)$, and
the connecting homomorphism $\der_{k,j}$ 
associated to the short exact sequence 
$	0\arr\Nd{BP_*/I_{k-1}} \arT2{v_{k-1}^j}BP_*/I_{k-1}\arr N_{k-1}(j)\arr 0$.
We abbreviate $\der_{k,1}$ to $\der_k$. 
Here we consider the Greek letter elements of $H^*BP_*/I_{n-1}$ defined by
\sk{
	\O \al^{(n-1)}_{t}\en u^t\in H^0BP_*/I_{n-1}\qand\\ 
	\al^{(n)}_{(t/j)}\en \der_{n,j}(v^t)\in H^1BP_*/I_{n-1}\mbx{for $v^t\in H^0N_{n-1}(j)$}
}
for $t>0$, and
\sk{
	\al_1\en \der_{1}(v_1)=h_0\in H^1BP_*\qand \be_1=\der_{1}\der_{2}(v_2)=b_0\in H^2BP_*.
}

\prop{al}{The elements $\al_1$ and $\be_1$ act on the Greek letter elements as follows:
\AR{l}{ \al_1\O \al_t^{(n-1)}\ne 0\in H^1BP_*/I_{n-1},\ak
	\be_1\O \al_t^{(n-1)}\ne 0\in H^2BP_*/I_{n-1};}
and if the Greek letter elements $\al^{(n)}_{(sp^i/j)}$ has an internal degree greater than $2(p^n-1)(e(n-1)-1)$, then
\AR{l}{ 
	\al_1\al^{(n)}_{(sp^i/j)}\ne 0\in H^2BP_*/I_{n-1}\ \mbox{if $[i]\ne 0$, $p\nmid (s+1)$ or $p^2\mid (s+1);$ and}\\
	\be_1\al^{(n)}_{(sp^i/j)}\ne 0\in H^3BP_*/I_{n-1}\ \mbox{if $n\ne 5$, $[i]\ne 1$ or $p\nmid (s+1)$.}\\
}
}

In order to prove this, we make a chromatic argument:
Let $N_k^0$ denote the $BP_*BP$-comodule $BP_*/I_k$,  
and put
$M_k^0=v_k^{-1}N_k^0$.
We denote the cokernel of the inclusion $N^0_k\to M^0_k$ by $N_k^1$, so that
$0\arr N^0_k\to \Nd{M^0_k} \ar{\ps} N_k^1\to 0$ is an exact sequence.
Let $\wt \der_{k+1}\cln H^sN_k^1\to H^{s+1}N_k^0$ be the connecting homomorphism associated to the short exact sequence.
We notice that $N^1_{k}={\rm colim}_j N_{k}(j)$ with inclusion $\ph_j\cln N_{k}(j)\to N^1_{k}$ given by $\ph_j(x)=x/u^j$, and that the connecting homomorphism
 $\der_{n,j}\cln H^sN_{n-1}(j)$ $\to H^{s+1}N_{n-1}^0$ factorizes to $\wt \der_n \ph_j$.

\lem{ap:0}{For an element $x_i^s/u^j\in H^0N^1_{n-1}$ for $0<j\le a_i$ $(j\le p^i$ if $s=1)$, $\al_1$ and $\be_1$ act on it as follows:
\AR{l}{
	x_i^s\al_1/u^j\ne 0  \in H^1N^1_{n-1}\mbx{ if $[i]\ne 0$, $p\nmid (s+1)$ or $p^2\mid (s+1)$; and}\\
	x_i^s\be_1/u^j\ne 0  \in H^2N^1_{n-1}\mbx{if $n\ne 5$, $[i]\ne 1$ or $p\nmid (s+1)$. }
} }

\pf{A change of rings theorem of Miller and Ravenel \cite{mr} shows that
the module $H^sM_{n-1}^1$ is isomorphic to $\e^sB$.
	By \kko{de:1}, we see that $x_i^sh_0/u\ne 0\in \e^1B$ unless $[i]=0$, $p\mid (s+1)$ and $p^2\nmid (s+1)$. This shows the first non-triviality.
Similarly,
since we have shown that  \kko{ses:ijs} is exact, we see that $x_i^s\be_1/u\ne 0\in \e^2B$ unless $n=5$, $[i]=1$ and $p\mid (s+1)$.
 }

\lem{ap:2}{Let $\xi_1$ denote $\al_1$ or $\be_1$, and $x\in H^0N_{n-1}^1$, and suppose that $x\xi_1$ has an  internal degree greater than $2(p^{n-1}-1)(e(n-1)-1)$. If $x\xi_1\in H^sN_{n-1}^1\ne 0$, then $\wt \der_n(x)\xi_1\ne 0\in H^{s+1}N_{n-1}^0$.
}

\pf{
It suffices to show that $x\xi_1$ is not in the image of $\ps_*\cln H^sM_{n-1}^0\to H^sN_{n-1}^1$.
Again the
 change of rings theorem shows that
the module $H^sM_{n-1}^0$ is isomorphic to the module of Lemma \ref{E1} with substituting $n-1$ for $n$. Note that every generator of it except for $\ze_{n-1}$ belongs to $H^sN_{n-1}^0$, and also is $u^{e(n-1)}\ze_{n-1}$ (\cf \cite{r:book}).
It follows that
every element of the image of $\ps_*$ has an internal degree no greater than  $2(e(n-1)-1)(p^{n-1}-1)$.
Thus the lemma follows.
}

\pfc{Proposition \ref{al}}{
The module $H^sM_{n-1}^0$ 
contains a submodule $k_*\lrk{h_0}$ if $s=1$ and $k_*\lrk{b_0}$ if $s=2$. 
Therefore, the first two relations hold.
The other relations follow from Lemmas \ref{ap:0} and \ref{ap:2}.
}

\pfc{Theorem \ref{1}}{
Note that
$\overline\al_t^{(3)}=\overline\gamma_t=v_3^t$, 
and we obtain the theorem  from Proposition \ref{al} at $n=4$.
}



\begin{thebibliography}{99}
\bibitem{as}
 Y. Arita and K. Shimomura, The chromatic $E_1$-term $H^1M^1_1$ at the prime $3$. Hiroshima Math. J. 26 (1996), 415--431.
\bibitem{h}
H.-W. Henn,
Centralizers of elementary abelian $p$-subgroups and mod-$p$ cohomology of profinite groups,
Duke Math. J. {\bf 91} (1998), 561--585.

\bibitem{hs}
H.~Hirata and K.~Shimomura, The chromatic $E_1$-term $H^1M^1_2$ for an odd prime, in preparation.

\bibitem{k}
D. Kraines,
	Massey higher products, Trans. Amer. Math. Soc. {\bf 124} (1966), 431--449.
\bibitem{may}
J. P. May,
	Matric Massey products, J. Alg. {\bf 12} (1969), 533--568.
\bibitem{may:st}
J. P. May,
A general algebraic approach to Steenrod operations, The Steenrod Algebra and its Applications, Lecture Notes in Mathematics {\bf 168} (1970), 153--231.
\bibitem{mr}
 H.~R.~Miller and D.~C.~Ravenel,
Morava stabilizer algebras and the localization of Novikov's $E_2$-term, Duke Math. J. , {\bf 44} (1977), 433-447. 
\bibitem{mrw}
 H.~R.~Miller, D.~C.~Ravenel, and W.~S.~Wilson,
     Periodic phenomena in Adams-Novikov spectral sequence,  Ann. of Math.           {\bf 106} (1977), 469--516.
\bibitem{n:n}
H.~Nakai, The chromatic $E_1$-term $H^0M^2_1$ for $p>3$, New York J. Math. {\bf 6} (2000), 21--54 (electronic). 
\bibitem{n:k}
H.~Nakai, The structure of ${\rm Ext}_{{\rm BP}_*{\rm BP}}^0({\rm BP}_*,M_1^2)$ for $p=3$, Mem. Fac. Sci. Kochi Univ. Ser. A Math. 23 (2002), 27-- 44.
\bibitem{r:M}
 D.~C.~Ravenel,
	The cohomology of the Morava stabilizer algebras, Math. Z. {\bf 152} (1977), 287--297.	
\bibitem{r:loc}
 D.~C.~Ravenel,	Localization with respect to certain periodic homology theories, Amer. J. Math., {\bf 106} (1984), 351--414.


\bibitem{r:orange}
D. C. Ravenel, Nilpotence and periodicity in stable homotopy theory, Annals of Mathematics
Studies, vol. 128, Princeton University Press, 1992.

\bibitem{r:book}
 D.~C.~Ravenel,
           {\it Complex cobordism and stable homotopy groups of spheres},                               AMS Chelsea Publishing, Providence, 2004.
\bibitem{s:86}
K. Shimomura,
On the Adams-Novikov spectral sequence and products of $\beta$-elements, Hiroshima Math. J. {\bf 16} (1986), 209--224. 
\bibitem{s-tottori}
K. Shimomura,
The chromatic $E_1$-term $H^1M^1_2$ and its application to the homology of the Toda-Smith spectrum $V(1)$, J. Fac. Educ. Tottori Univ. (Nat. Sci.) 39 (1990), 63--83. Corrections to ``The chromatic $E_1$-term $H^1M^1_2$ and its application to the homology of the Toda-Smith spectrum $V(1)$'', J. Fac. Educ. Tottori Univ. (Nat. Sci.) 41 (1992), 7--11. 
\bibitem{s:90}
K. Shimomura,
The chromatic $E_1$-term $H^0M^2_n$ for $n> 1$, J. Fac. Educ. Tottori Univ. (Nat. Sci.) {\bf 39} (1990), 103--121
\bibitem{s:97}
K. Shimomura,
The homotopy groups of $L_2$-localized Toda-Smith spectrum $V(1)$ at the prime 3, Trans.~Amer.~Math.~Soc. {\bf 349} (1997), 1821--1850. 
\bibitem{s:00}
K. Shimomura,
The homotopy groups of the $L_2$-localized mod 3 Moore spectrum, J. Math. Soc. of Japan {\bf 51} (2000), 65--90.
\bibitem{st}
K. Shimomura and H. Tamura, 
Non-triviality of some compositions of $\beta$-elements in the stable homotopy of the Moore spaces, Hiroshima Math. J. {\bf 16} (1986), 121--133. 
\bibitem{sy}
K. Shimomura and A. Yabe,
The homotopy groups $\pi_*(L_2S^0)$, Topology {\bf 34} (1995), 261--289. 
\bibitem{sw}
K. Shimomura and X. Wang,
The homotopy groups $\pi_*(L_2S^0)$ at the prime 3 , Topology {\bf 41} (2002), 1183--1198. 
\bibitem{t}
H. Toda, On spectra realizing exterior parts of the Steenrod algebra, Topology 10 (1971), 53-65.
\end{thebibliography}
\end{document}